\magnification1200

\centerline{\bf Trees of manifolds as boundaries of spaces and groups}

\medskip
\centerline{{ Jacek \'Swi\c atkowski}
\footnote{*}{This work was supported by the Polish Ministry 
of Science and Higher Education (MNiSW), grant N201 541738,
and by the Polish National Science Centre (NCN), grant 2012/06/A/ST1/00259.}}

% 19 stycznia 2016

\bigskip
{\bf Abstract.} We show that trees of manifolds, the topological spaces
introduced by Jakobsche, appear as boundaries at infinity of various
spaces and groups. In particular, they appear as Gromov boundaries
of some hyperbolic groups, of arbitrary dimension, obtained
by the procedure of strict hyperbolization.
We also recognize these spaces as boundaries of arbitrary Coxeter groups 
with manifold nerves, and as Gromov boundaries of the fundamental
groups of singular spaces obtained from some finite volume hyperbolic
manifolds by cutting off their cusps and collapsing the resulting
boundary tori to points.


\bigskip
\centerline{\bf Introduction}

\medskip
Not many explicit topological spaces are known to be homeomorphic
to Gromov boundaries of hyperbolic groups. The list consists of spheres
and Sierpi\'nski compacta of arbitrary dimension, Menger compacta
in dimensions $d\in\{ 1,2,3 \}$, Pontriagin sphere and Pontriagin surfaces
(which are 2-dimensional), and some trees of 3-manifolds (which are
certain spaces of
topological dimension 3). See [PS] for more detailed comments concerning
this list. 
One of the goals of this paper is to extend this list to trees
of $n$-manifolds in arbitrary dimension $n$ (for which the methods
used in [PS], in the case $n=3$, are insufficient).

Trees of manifolds have been formally introduced by W\l odzimierz
Jakobsche in [J], though their idea goes back to Ancel and Siebenmann [AS].
These are certain explicit homogeneous metric compacta, typically not ANR,
appearing in abundance in arbitrary finite topological dimension.
We describe them in detail in Section 1, and here we only mention that
each closed connected topological $n$-manifold $M$ determines 
uniquely one such space,
denoted ${\cal X}(M)$, which has topological dimension $n$, and which
we call {\it the tree of manifolds $M$.}
Our first result is the following (compare Theorem 5.1 in the text,
which is slightly more general).

\medskip\noindent
{\bf Theorem 1.}
\item{(1)}
{\it For any natural number $n$, if a closed connected 
orientable PL $n$-manifold $M$ bounds
a compact orientable PL $(n+1)$-manifold, 
then the tree of manifolds ${\cal X}(M)$
is homeomorphic to the Gromov boundary of some hyperbolic group.}
\item{(2)}
{\it For any closed connected 
non-orientable PL manifold $N$ the tree of manifolds ${\cal X}(N)$
is homeomorphic to the Gromov boundary of some hyperbolic group.}

\medskip
Appropriate hyperbolic groups as in the assertion of Theorem 1
will be constructed in Section 5,
using the procedure of strict hyperbolization
due to Charney and Davis [CD].
Theorem 1 follows from a more general result, Main Theorem,
formulated in Section 2 (and proved in Sections 3 and 4), 
which concerns boundaries of certain $CAT(\kappa)$
pseudomanifolds, $\kappa\le0$, called $({\cal M},\kappa)$-pseudomanifolds. Another application of 
this theorem, discussed in detail in Section 5, 
is the following correction and extension 
of the result of H. Fischer [F] (who has studied only the case of
right-angled Coxeter groups with orientable manifold nerves).

\medskip\noindent
{\bf Theorem 2.}
{\it Let $(W,S)$ be a Coxeter system (not necessarily right angled)
whose nerve is a 
PL triangulation of a closed connected manifold $M$. Then the boundary
at infinity
of $(W,S)$ (i.e. the visual boundary of the corresponding Davis-Moussong
complex for $(W,S)$) is homeomorphic to the following tree of manifolds:}

\item{(1)} {\it ${\cal X}(M\#\overline M)$, if $M$ is orientable, where 
$M,\overline M$ are two oppositely oriented copies of $M$;}

\item{(2)} {\it ${\cal X}(M)$, if $M$ is non-orientable.}

\medskip
Proofs of Theorems 1 and 2 (or even slightly more general results) 
are presented in Section 5.
A variation of Main Theorem presented in Section 6 (Theorem 6.2)
yields the following result, which applies to a class of hyperbolic groups
(having representatives in arbitrary dimension) constructed 
by L. Mosher and M. Sageev in [MS] (see also [FM]).

\medskip\noindent
{\bf Theorem 3.}
{\it Let $M$ be a finite volume complete non-compact hyperbolic 
$(n+1)$-manifold with toral cusps. Suppose that after removing some
open horoball neighbourhoods of all cusps we get a compact
$(n+1)$-manifold $M^\circ$ whose euclidean toral boundary
components contain no closed geodesics of length $\le 2\pi$.
Let $\Gamma$ be a hyperbolic group which is the fundamental
group of a pseudomanifold obtained from $M^\circ$
by collapsing all its boundary components to points.
Then the Gromov boundary $\partial\Gamma$ is homeomorphic
to the tree of tori ${\cal X}(T^n)$.}

\medskip
The proofs of the results in this paper required some new ideas
and tools. One of them is a new and more flexible characterization 
of trees of manifolds as limits of inverse sequences of manifolds,
given in Definition 1.1, and established by the author 
in a separate paper [Sw]. Remarks after Definition 1.1 clarify
the novelty of this characterization.
Another new ingredient is a more careful or pushed further 
(than in [DJ], [F] and [FM])
analysis of geodesic projections between concentric spheres
in $CAT(\kappa)$ pseudomanifolds. 
This analysis allows to approximate projections as above by more regular 
maps, which in turn allows to understand the boundaries of the
corresponding pseudomanifolds.

The above mentioned analysis of geodesic projections builds upon 
some results of M. Davis and T. Januszkiewicz in [DJ]. Since the proofs
of these results given in [DJ] are incomplete, as we were informed by 
the referee of the present paper, we include an appendix. This appendix
presents complete proofs of these results, in the form provided later
by the authors of [DJ], but never published. The authors of [DJ] have
approved such an arrangement.

\smallskip
The author thanks Mike Davis, Tadeusz Januszkiewicz and Damian Osajda
for helpful and inspiring discussions, and the referee for suggestions of
significant improvements of the paper.

\bigskip\noindent
{\bf 1. Trees of manifolds.}

\medskip
Trees of manifolds have been introduced by W\l odzimierz Jakobsche in [J].
These are some homogeneous metric compacta, typically not ANR, appearing
in abundance in any finite topological dimension. We briefly recall
their description, in the original setting of Jakobsche modified
and extended as in [Sw]. Although the description is valid for arbitrary
topological manifolds, in this paper we make use only of the subclass
related to PL manifolds.

\medskip\noindent
{\bf Definiton 1.1.}
Let $\cal M$ be a finite family of closed connected $n$-dimensional 
topological manifolds,
either all oriented, or at least one of which is non-orientable.
Let $${\cal J}=(\{ X_i:i\ge1 \},\{ \pi_i:i\ge1 \})$$
be an inverse sequence consisting of closed connected topological
$n$-manifolds $X_i$
and maps $\pi_i:X_{i+1}\to X_i$. Assume furthermore that if all the manifolds in
$\cal M$ are oriented then all $X_i$ are also oriented. We say that $\cal J$ is a
{\it weak Jakobsche inverse sequence for $\cal M$} if 
for all $i\ge1$ and all $M\in{\cal M}$ one can choose finite families
${\cal D}_i$ of collared $n$-disks in $X_i$, partitioned into subfamilies
${\cal D}_i^M:M\in{\cal M}$, such that:

\smallskip
\itemitem{(1)} for each $i\ge1$ the disks in the family ${\cal D}_i$
are pairwise disjoint;

\itemitem{(2)} for each $i\ge1$ the map $\pi_i$ maps the preimage
$\pi_i^{-1}(X_i\setminus\cup\{ \hbox{int}(D):D\in{\cal D}_i \})$
homeomorphically onto $X_i\setminus\cup\{ \hbox{int}(D):D\in{\cal D}_i \}$;

\itemitem{(3a)} $X_1$ is homeomorphic to one of the manifolds from $\cal M$,
and if the manifolds in $\cal M$ are oriented, we require that this 
homeomorphism respects orientations;

\itemitem{(3b)} for each $i\ge1$, for each $M\in{\cal M}$, and
for any $D\in{\cal D}_i^M$
the preimage $\pi_i^{-1}(D)$ is homeomorphic to $M\setminus\hbox{int}(\Delta)$,
where $\Delta$ is some collared $n$-disk in $M$; furthermore, if the manifolds
in $\cal M$ are oriented, we require that the above homeomorphism
respects the orientations induced from $X_{i+1}$ and from $M$;

\itemitem{(4)} if $i<j$, $D\in{\cal D}_i$, $D'\in{\cal D}_j$, then
$\pi_{i,j}(D')\cap\partial D=\emptyset$, where
$\pi_{i,j}:=\pi_{i}\circ\pi_{i+1}\circ\dots\circ\pi_{j-1}$;

\itemitem{(5)} for each $i\ge1$ the family 
$\{ \pi_{i,j}(D):j\ge i, D\in{\cal D}_j \}$
of subsets of $X_i$ is {\it null}, i.e. the diameters of these subsets
converge to 0; here $\pi_{i,i}$ denotes the identity map
on $X_i$;

\itemitem{(6)} for any $i\ge1$ and each $M\in{\cal M}$ the set
$\bigcup_{j=i}^\infty\pi_{i,j}(\cup{\cal D}_j^M)$
is dense in $X_i$.

\medskip\noindent
{\bf Remarks.}

\item{(1)} It follows from conditions (1), (2), (3a) and (3b)
that each $X_i$ is the connected sum of a family of manifolds
each homeomorphic to one of the manifolds in $\cal M$;
moreover, if the manifolds in $\cal M$ are oriented, the above
mentioned homeomorphisms (and the involved operation 
of connected sum) respect the orientations.

\item{(2)} In the case when the manifolds in $\cal M$ are oriented,
conditions (1)-(5) in Definition 1.1 basically coincide with conditions
(1)-(6) in [J], Section 2, p. 82.

\item{(3)} Condition (6) in Definition 1.1 implies condition (7)
in [J], but it is essentially weaker than the conjunction of conditions 
(7) and (8) of [J] (except 
when the family $\cal M$ consists of a single manifold $M$,
in which case condition (6) of Definition 1.1 is equivalent to
the conjunction of conditions (7) and (8) of [J], as it was
observed and exploited in [F] and [Z]).

\item{(4)} Definition 1.1 describes exactly the class of inverse
sequences  naturally associated
to weakly saturated tree systems of manifolds from $\cal M$,
as described in [Sw], Section 3.E. It represents significant, and probably
close to optimal,
relaxation of the initial set of conditions provided
by Jakobsche in [J], for which the inverse limit of the corresponding sequence
depends uniquely on $\cal M$ (see Theoerem 1.3 below).

\medskip
The following result is a reformulation of Corollary 3.E.4 of [Sw].

\medskip\noindent
{\bf Theorem 1.3.}
{\it Given $\cal M$ as in Definition 1.1, any two weak
Jakobsche inverse sequences for $\cal M$ 
have homeomorphic inverse limits.}

\medskip\noindent
{\bf Definition 1.4.}
Given $\cal M$ as in Definition 1.1, denote by ${\cal X}(\cal M)$,
and call {\it the tree of manifolds from $\cal M$}, the space 
homeomorphic to the inverse limit of some (and hence any)
weak Jakobsche inverse sequence for $\cal M$.

\medskip\noindent
{\bf Remarks.}

\item{(1)} If ${\cal M}=\{ M \}$, we denote the corresponding space
${\cal X}({\cal M})$ simply by ${\cal X}(M)$, and call it 
{\it the tree of manifolds $M$}.

\item{(2)} If ${\cal M}=\{ M_1,\dots,M_k \}$, it is known that the space
${\cal X}({\cal M})$ is homeomorphic to the tree of manifolds
$M_1\#\dots \#M_k$, i.e. to the space ${\cal X}(M_1\#\dots \#M_k)$
(see e.g. Corollary 3.E.4 in [Sw]).

\medskip
As it was shown in [J] and [St], trees of manifolds 
${\cal X}({\cal M})$ are connected homogeneous metric compacta of
topological dimension equal to the dimension of the manifolds in $\cal M$.
As it was observed in Corollary 3.3 of [PS], trees of manifolds
in dimensions $\ge2$ have no local cut points, and in fact they are
{\it Cantor manifolds}, i.e. no subsets of topological codimension
$\ge2$ separate them.

Trees of manifolds ${\cal X}({\cal M})$ of the same topological
dimension $n$ can be sometimes distinguished by means of
their homotopical or homological invariants, or the shape theoretic
invariants. For example, it is known that if $n\ge3$ 
then the shape fundamental group
$\check{\pi}_1({\cal X}({\cal M}))$ is isomorphic to the inverse limit
of the increasing free products $G_1*G_2*\dots*G_k$, where the
additional factors in larger products are collapsed to identity
while the common factors are mapped to each other through identities,
and where groups $G_i$ are the copies of the fundamental groups 
$\pi_1(M):M\in{\cal M}$, each group appearing infinitely often
(see [FG]). Obviously, this shape fundamental group 
is sometimes sufficient to distinguish
some trees of manifolds. The general question of classifying trees
of manifolds up to homeomorphism remains open.

\bigskip

\noindent
{\bf 2. $({\cal M},\kappa)$-pseudomanifolds and Main Theorem.}

\medskip
Let $\cal M$ be a finite collection of (pairwise distinct) 
closed connected PL manifolds
of the same dimension $n$, distinct from the standard PL $n$-sphere,
and let $\kappa\in\{ 0,-1 \}$.
An $({\cal M},\kappa)$-{\it pseudomanifold}
is a metric polyhedral complex $X$ of 
piecewise constant curvature $\kappa$,
with finite shapes, and such that

\item{(1)} $X$ is $CAT(\kappa)$;

\item{(2)} for each point $x\in X$ the link $\hbox{Lk}(x,X)$
(viewed as combinatorial polyhedral complex)
is either a PL $n$-sphere, or is PL-homeomorphic to some $M\in{\cal M}$;
furthermore, if all $M\in{\cal M}$ are oriented, we assume that $X$
is also oriented, and that each homeomorphism $\hbox{Lk}(x,X)\to M$
as above respects the orientations;

\item{(3)} for each $M\in{\cal M}$ the set 
$\Lambda_M:=\{ x\in X:\hbox{Lk}(x,X)\cong M \}$
is a net in $X$, i.e. for some $R>0$ each ball of radius $R$ in $X$
intersects $\Lambda_M$ 
(if all $M\in{\cal M}$ are oriented, the symbol $\cong$ above
denotes relation of being PL homeomorphic as oriented manifolds). 

\medskip

Note that,
since $X$ above has finite shapes, it is a geodesic metric space,
so that it makes sense to speak of the $CAT(\kappa)$ property for it.
Moreover,
by condition (2), $X$ is automatically
an $(n+1)$-dimensional pseudomanifold. 
The assumption that the link of every point of $X$ is a manifold implies
that only the vertices of $X$ can have links that are not spheres.
Hence, if $x\in X$ and $\hbox{Lk}(x,X)\cong M\in{\cal M}$,
than $x$ must be a vertex of $X$.
The finite shapes property implies that the vertices of $X$ form a discrete
set. Therefore, each set $\Lambda_M$
as in condition (3) above is discrete.
Putting $\Lambda:=\cup_{M\in{\cal M}}\Lambda_M$, we get
that $\Lambda$, which will be called {\it the singular set} of $X$, 
is also discrete.

We are now ready to state Main Theorem, which is the main result
of the paper. Its proof occupies Sections 3 and 4. In Section 5
we show few classes of examples to which this theorem applies,
in particular getting the proofs of Theorems 1 and 2 of the introduction.

\bigskip
\noindent
{\bf Main Theorem.}
{\it The visual boundary $\partial X$ of any 
$({\cal M},\kappa)$-pseudomanifold $X$ is homeomorphic
to the tree of manifolds ${\cal X}({\cal M})$.}

\bigskip
\noindent
{\bf 3. Proof of Main Theorem.}

\medskip

Fix any point $x_0\in\Lambda$.
Denote the spheres and the closed and open balls of radii $R$ in $X$,
centered at $x_0$, respectively, by $S_R=\{ x\in X:d(x,x_0)=R \}$,
$B_R=\{ x\in X:d(x,x_0)\le R \}$ and 
$B^\bullet_R=\{ x\in X:d(x,x_0)<R \}$.
By discreteness of $\Lambda$, the set of numbers
$\{ d(x_0,p):p\in\Lambda \}$ is also discrete. Order this set
into the increasing sequence $R_i:i\ge0$, with $R_0=0$.
For each $i\ge0$ let $p_{i,1},\dots,p_{i,k_i}$ be the points
of the intersection $S_{R_i}\cap\Lambda$.
(In particular, we have $k_0=1$ and $p_{0,1}=x_0$.)

For each $i\ge1$, let $G_i:X\setminus B^\bullet_{R_i}\to S_{R_i}$
be the geodesic projection towards $x_0$, i.e. the map which
to any point $x\in X\setminus B^\bullet_{R_i}$ associates
the unique point $x'$ in the intersection of the sphere $S_{R_i}$
with the geodesic segment $[x,x_0]$. Denote also by $g_i:S_{R_{i+1}}\to S_{R_i}$
the restriction of $G_i$ to $S_{R_{i+1}}$.
Consider the inverse sequence ${\cal S}=(\{ S_{R_i} \},\{ g_i \})$
and recall that $$\partial X=\lim_{\longleftarrow}{\cal S}$$
(see II.8.5 in [BH], or the comments after Definition (2b.1) in [DJ]). 
Before getting to the core of the proof,
we need few preparatory results.

\medskip
The next result was essentially proved in Section 3 of [DJ].
Its special case, in which $X$ is an $({\cal M},0)$-pseudomanifold
and $\cal M$ consists of homology spheres, appears as
Theorem (3d.1)(i) in [DJ]. The proof of this theorem in [DJ]
relies on the argument given in the paragraph starting at the bottom
of page 371 and ending at page 372 of that paper; this argument 
is incomplete (as we were informed by the referee of the present paper),
and thus we provide its complete version in the appendix to the present
paper. In particular, the lemma below is a special case of Lemma B
of the appendix.

\medskip\noindent
{\bf Lemma 3.1.}
{\it Let $X$ be an $(n+1)$-dimensional $({\cal M},\kappa)$-pseudomanifold.
Then each sphere $S_R(x,X)$ in $X$ (with positive radius $R$)
is a closed $n$-manifold.}

\medskip
The next result is a special case of Lemma C of the appendix to the
present paper. The notion of a cell-like map appearing in the statement
of this result is recalled in the appendix (right after the statement
of Lemma C).

\medskip\noindent
{\bf Lemma 3.2.}
{\it Given an $({\cal M},\kappa)$-pseudomanifold $X$, any $R'>R>0$,
and any $x\in X$,
let $g:S_{R'}(x,X)\to S_R(x,X)$ be the geodesic projection
towards $x$.
If the set $B^\bullet_{R'}(x,X)\setminus B^\bullet_R(x,X)$
contains no singular point from $\Lambda$ then $g$ is a cell-like map.
In particular, in the notation established at the beginning of this
section, for any $0<\epsilon<R_{i+1}-R_i$ the geodesic projection
$\phi_i:S_{R_{i+1}}\to S_{R_i+\epsilon}$ towards $x_0$ is cell-like.}

\medskip
Recall that, by a result of M. Brown [Br]
(to which we refer as {\it Brown's Lemma}), if we replace the maps $g_i$
in the inverse sequence $\cal S$, successively, by sufficiently close
maps $g_i':S_{R_{i+1}}\to S_{R_i}$ then the limit of the resulting
inverse sequence ${\cal S}'=(\{ S_{R_i} \},\{ g_i' \})$ remains unchanged
(up to homeomorphism),
i.e. $\lim_{\longleftarrow}{\cal S}'\cong \lim_{\longleftarrow}{\cal S}$.
The next proposition describes particularly nice approximations $g_i'$, 
as above, of the maps $g_i$.

\medskip\noindent
{\bf Proposition 3.3.} {\it Let $\dim(X)=n+1$.
For all $i\ge1$ there exist maps $g_i':S_{R_{i+1}}\to S_{R_i}$ 
satisfying the following: }

\smallskip
\item{(1)} {\it $g_i'$ are so close to $g_i$ that, according to Brown's Lemma,
we have 
$$\lim_{\longleftarrow}{\cal S}'\cong \lim_{\longleftarrow}{\cal S}
\cong\partial X;$$} 

\item{(2)} {\it for each $i\ge1$, putting 
$S_{R_i}^\circ=S_{R_i}\setminus\{ p_{i,1},\dots,p_{i,k_i}\}$,
the restriction 
$$g_i'|_{(g_i')^{-1}(S_{R_i}^\circ)}:(g_i')^{-1}(S_{R_i}^\circ)
\to S_{R_i}^\circ$$
is a homeomorphism;}

\item{(3)} {\it for each $i\ge1$ and each $1\le j\le k_i$, if $D_{i,j}$
is a sufficiently small collared $n$-disk in $S_{R_i}$
such that $p_{i,j}\in\hbox{int}(D_{i,j})$
then $$(g_i')^{-1}(D_{i,j})\cong M\setminus\hbox{int}(\Delta),$$
where $M=\hbox{Lk}(p_{i,j},X)\in{\cal M}$, and where $\Delta$ is any
collared $n$-disk in $M$;}

\item{(4)} {\it for each $M\in{\cal M}$ and for all $i\ge1$,
the images in $S_{R_i}$ of the points from the set 
$\Lambda_M\cap(X\setminus B_{R_i})$,
through appropriate compositions of the maps $g_k':k\ge i$,
form a dense subset of $S_{R_i}$.}

\medskip
We postpone the proof of Proposition 3.3 until the next section,
and in the remaining part of this section we 
complete the proof of Main Theorem, using the proposition.
We do this by showing that an inverse sequence 
${\cal S}'=(\{ S_{R_i} \},\{ g_i' \})$ resulting from
Proposition 3.3 is a weak  Jakobsche inverse sequence for $\cal M$.
More precisely, we describe families ${\cal D}_i^M$ of disks
as required in Definition 1.1, inductively with respect to $i$.

We start with checking condition (3a) of Definition 1.1, i.e. 
showing that $S_{R_1}$ is homeomorphic to one of the manifolds 
from $\cal M$. Here, and later in Section 4, we will need
the concept of a {\it cone neighbourhood} of a point
in a piecewise constant curvature polyhedral complex.
Let $Y$ be a metric polyhedral complex of constant curvature
$k$ equal $1$, $-1$ or $0$, with finite shapes. 
Then for any point $y\in Y$,
and any sufficiently small $\epsilon>0$,
the ball $B_\epsilon(y,Y)$ isometrically coincides with the ball of the same 
radius $\epsilon$ in the $k$-cone $C_k(\hbox{Lk}(y,Y))$
centered at the cone vertex (see 
Definition I.5.6 on p. 59 of [BH]
or [Da], p. 505 for the 
definition of the $k$-cone, and
Theorem I.7.16 on p. 103 of [BH] 
for justification of the above claim). A {\it cone neighbourhood}
of $y$ in $Y$ is any ball $B_\epsilon(y,Y)$ as above.
An obvious (and useful) property 
is that geodesics started at any point of $Y$ do not bifurcate inside
a cone neighbourhood of this point. It follows that the natural map
from the boundary sphere $S_\epsilon(y,Y)$ of any cone neighbourhood,
to the link $\hbox{Lk}(y,Y)$, is a homeomorphism, which we will
view as the natural identification of the sphere with the link.

Coming back to checking condition (3a), choose $\epsilon$ so small
that $B_\epsilon$ is a cone neighbourhood of $x_0$ in $X$.
Then the sphere $S_\epsilon$ is homeomorphic to the link 
$\hbox{Lk}(x_0,X)$, i.e. to some manifold $M_0\in{\cal M}$.
By Lemma 3.2, the geodesic projection $\phi_0:S_{R_1}\to S_\epsilon$
is a cell-like map. By Approximation Theorem concerning cell-like maps
(recalled in the appendix), $\phi_0$ can be approximated by
homeomorphisms. It follows that  
$S_{R_1}$ is homeomorphic to $S_\epsilon$,
and hence also to $M_0$, as required.

We now turn to describing families ${\cal D}_i^M$ of disks
as required in Definition 1.1.
Fix an auxilliary sequence $\epsilon_i$ of positive real numbers,
converging to 0.
To start the inductive description, choose a family of pairwise disjoint
collared disks $D_{1,j}:1\le j\le k_1$ in $S_{R_1}$  such that
for each $j$ the following conditions hold:

\smallskip

\itemitem{(d1)} $D_{1,j}$ contains the point $p_{1,j}$ in its interior; 

\itemitem{(d2)} the diameter of $D_{1,j}$
is less than $\epsilon_1$;

\itemitem{(d3)} 
the preimage $(g_1')^{-1}(D_{1,j})$ is homeomorphic to
$M\setminus\hbox{int}(\Delta)$, 
where $M\in{\cal M}$ is homeomorphic to the link $\hbox{Lk}(p_{1,j},X)$
and $\Delta$ is a
collared $n$-disk in $M$; 

\itemitem{(d4)} denoting by $\Lambda(1)$ the set of images in $S_{R_1}$,
through appropriate compositions of the maps $g_l':l\ge1$,
of all points in $\Lambda\cap(X\setminus B_{R_1})$, we have that
$\Lambda(1)\cap\partial D_{1,j}=\emptyset$.

\smallskip
\noindent
A choice as above is possible due to property (3) in Proposition 3.3
(which allows to obtain (d3)),
and due to the fact that the set $\Lambda(1)$ is countable
(which allows to get (d4)).
Put ${\cal D}_1=\{ D_{1,j}:1\le j\le k_1 \}$, and split this set
into subsets ${\cal D}_1^M$ according to the rule that
$D_{1,j}\in{\cal D}_1^{M}$ iff ${Lk}(p_{1,j},X)\cong M$.

Now, suppose that $m\ge2$ and that for all $1\le i< m$,
and all $M\in{\cal M}$, we have already
described the families ${\cal D}_i^M$
satisfying conditions (1), (2), (3b) and (4) of Definition 1.1.
Suppose also that any disk in a so
far described family ${\cal D}_i^M$,
as well as its image in any $S_{R_k}:k<i$ through
the appropriate composition of the maps $g_l'$,
has diameter less than $\epsilon_i$.
(We demand this to ensure that condition (5) of Definition 1.1
holds, after describing all families of disks.)
Clearly, if $m=2$ then these assumptions hold, which can be
easily deduced from the first step of the construction
given in the previous paragraph.

Choose a family of pairwise disjoint
collared disks $D_{m,j}:1\le j\le k_m$ in $S_{R_m}$  such that
for each $j$ the following conditions hold:

\smallskip

\itemitem{(d1')} $D_{m,j}$ contains the point $p_{m,j}$ in its interior; 

\itemitem{(d2')} the diameter of $D_{m,j}$, 
as well as of its image in any $S_{R_k}:k<m$ through the appropriate
composition of the maps $g_l'$,
is less than $\epsilon_m$;

\itemitem{(d3')} 
the preimage $(g_m')^{-1}(D_{m,j})$ is homeomorphic to
$M\setminus\hbox{int}(\Delta)$, 
where $M\in{\cal M}$ is homeomorphic to the link $\hbox{Lk}(p_{m,j},X)$
and $\Delta$ is a
collared $n$-disk in $M$; 

\itemitem{(d4')} denoting by $\Lambda(m)$ the set of images in $S_{R_m}$,
through appropriate compositions of the maps $g_l':l\ge m$,
of all points in $\Lambda\cap(X\setminus B_{R_m})$, we have that
$\Lambda(m)\cap\partial D_{m,j}=\emptyset$;

\itemitem{(d5')} for any $k<m$ and any $D\in{\cal D}_k$
the image of $D_{m,j}$ in $S_{R_k}$, through the appropriate
composition of the maps $g_l'$, is disjoint from $\partial D$.

\smallskip
\noindent
The fact that choosing $D_{m,j}$ sufficiently small we can fulfill
(d5') follows from having condition (d4') satisfied,
by inductive hypothesis,
for $m$ replaced with any $i<m$.
(Clearly, condition (d5') corresponds to condition (4) in Definition 1.1.)

The families of disks obtained by the above inductive description
clearly satisfy conditions (1)-(5) of Definition 1.1. It remains
to check condition (6) of this definition. However, this condition
follows fairly directly from condition (4) in Proposition 3.3
(which holds for the maps $g_i'$ in the inverse sequence ${\cal S}'$),
and from the fact that each point $p_{i,j}\in\Lambda_M$ is included
in the disk $D_{i,j}\in{\cal D}_i^M$.

This completes the proof of Main Theorem.

\bigskip\noindent
{\bf 4. Proof of Proposition 3.3.}

\medskip
The proof of Proposition 3.3 is split into a series of auxilliary
observations and partial results, and it is completed in the last part
of the section.
We use the notation established 
at the beginning of Section 3,
where $X$ is an $({\cal M},\kappa)$-pseudomanifold, $\Lambda$ is
its singular set, $S_{R_i}$ are the spheres in $X$ centered at a fixed
point $x_0\in\Lambda$, of appropriately chosen radii $R_i$,
and $G_i:X\setminus B^\bullet_{R_i}\to S_{R_i}$ are the geodesic
projections towards $x_0$.  

The following observation is an easy consequence of the assumptions
that $X$ is $CAT(\kappa)$ for some $\kappa\in\{0,-1\}$, and
that for each $M\in{\cal M}$ the singular set $\Lambda_M$ is a net
in $X$ (we omit the proof).

\medskip\noindent
{\bf Claim 4.1.}
{\it For each $M\in{\cal M}$ and for each $i\ge1$ the image set
$G_i(\Lambda_M\cap(X\setminus B_{R_i}))$ is dense in $S_{R_i}$.}

\medskip 
The next result follows essentially by the same arguments as those used in
the proof of Lemma C in the appendix.
More precisely, the same arguments show that point preimages
of the restricted map $\psi_i$ are cell-like sets. Moreover, the fact that
$\psi_i$ is a surjective map between compact spaces easily implies
that its restriction of the form as in the statement below is surjective and proper. We omit further details of the justification, but we make an additional
important observation that, in view of Lemma 3.1, the sets
$S_{R_i}\setminus\Lambda_{R_i}$ and  
$\psi_i^{-1}(S_{R_i}\setminus\Lambda_{R_i})$ appearing in the
statement below are open subsets in manifolds, and hence are
themselves manifolds.

\medskip\noindent
{\bf Lemma 4.2.}
{\it Under notation established at the beginning of Section 3,
put $\Lambda_{R_i}:=\Lambda\cap S_{R_i}$.
Given some $\epsilon\in(0,R_{i+1}-R_i)$, denote by
$\psi_i:S_{R_i+\epsilon}\to S_{R_i}$ the geodesic projection
towards $x_0$, i.e. the appropriate restriction of the map $G_i$.
Then the restriction of $\psi_i$ to the map from
$\psi_i^{-1}(S_{R_i}\setminus\Lambda_{R_i})$ to $S_{R_i}\setminus\Lambda_{R_i}$ is a cell-like map between manifolds.}

\medskip
Next three lemmas deal with approximations of geodesic projections
by the maps of better properties. They prepare the ground for 
the construction of approximations $g_i'$ as required in Proposition 3.3.
More precisely, Lemma 4.3 is related to condition (2) of the proposition,
while Lemmas 4.4 and 4.5 concern conditions (3) and (4), respectively.

\medskip\noindent
{\bf Lemma 4.3.}
{\it Let $\epsilon\in(0,R_{i+1}-R_i)$, and let 
$\psi_i:S_{R_i+\epsilon}\to S_{R_i}$ be the geodesic projection
towards $x_0$ in $X$, i.e. the appropriate restriction of the
map $G_i$. Then $\psi_i$ can be approximated by the maps
$\psi_i'$ such that:}

\smallskip
\item{(1)} {\it for each $1\le j\le k_i$ we have 
$(\psi_i')^{-1}(p_{i,j})=\psi_i^{-1}(p_{i,j})$;}

\item{(2)} {\it $\psi_i'$ restricted to the preimage
$(\psi_i')^{-1}(S_{R_i}\setminus\Lambda_{R_i})$
maps this set homeomorphically onto $S_{R_i}\setminus\Lambda_{R_i}$.}

\medskip\noindent
{\bf Proof:} Let 
$\bar\delta:S_{R_i+\epsilon}\to[0,\infty)$
be a continuous function such that 
$\bar\delta^{-1}(0)=\psi_i^{-1}(\Lambda_{R_i})$,
and let $\delta:S_{R_i+\epsilon}\setminus\psi_i^{-1}(\Lambda_{R_i})\to R_+$ 
be the restriction of $\bar\delta$.
Since a cell-like map between manifolds is a near-homeomorphism
(see Approximation Theorem in the appendix), 
it follows from Lemma 4.2 that there is a homeomorphism 
$\psi':S_{R_i+\epsilon}\setminus\psi_i^{-1}(\Lambda_{R_i})\to
S_{R_i}\setminus\Lambda_{R_i}$ such that
$d(\psi'(y),\psi_i(y))\le\delta(y)$ for all $y$ in the domain
of $\psi'$. By this estimate, and since $\delta(y)$ converges to $0$
as $y$ converges to a point in $\psi_i^{-1}(\Lambda_{R_i})$, 
$\psi'$ can be extended to a continuous map 
$\psi_i':S_{R_i+\epsilon}\to S_{R_i}$, by puting $\psi_i'(z)=\psi_i(z)$
for all $z\in \psi_i^{-1}(\Lambda_{R_i})$. 
This extension necessarily satisfies properties (1) and (2),
and clearly we can get in this way a map as close to $\psi_i$
as we wish. This finishes the proof.

\medskip\noindent
{\bf Lemma 4.4.}
{\it Under notation of Lemma 4.3, if $\epsilon$ is sufficiently small
then any map $\psi_i'$ fulfiling assertions (1) and (2) of this lemma
satisfies also the following:}

\smallskip
\item{($\star$)} {\it for each $1\le j\le k_i$ and any sufficiently small
collared $n$-disk $D$ in $S_{R_i}$ containing $p_{i,j}$ in its interior,
the preimage $(\psi_i')^{-1}(D)$ is homeomorphic to
$M\setminus\hbox{\rm int}(\Delta)$, where 
$M=\hbox{\rm Lk}(p_{i,j},X)\in{\cal M}$, and
where $\Delta$ is a collared $n$-disk in $M$.}

\medskip\noindent
{\bf Proof:}
Let $\epsilon$ be so small that for each $1\le j\le k_i$ the preimage
$C_{j}:=\psi_i^{-1}(p_{i,j})$, viewed as a subset of $X$, 
is contained in 
the interior of some cone neighbourhood
$B_{r_j}(p_{i,j},X)$. Consider the canonical identification of each
sphere $L_{j}:=S_{r_j}(p_{i,j},X)$ with the link 
$L_{j}^*:=\hbox{Lk}(p_{i,j},X)$.
Put $\{ w_j \}=[p_{i,j},x_0]\cap L_{j}$, and denote by $w_j^*$
the point in $L_j^*$ corresponding to $w_j$ under the above identification.
Let $A_j^*:=\{ y\in L_j^*:d_j(y,w_j^*)\ge\pi \}$, where $d_j$ is
the piecewise spherical metric in $L_j^*$, be the 
{\it infinitesimal shadow} of the point $p_{i,j}$ with respect to
$w_j^*$, and let $A_j$ be the corresponding subset in $L_j$.
Note that, by the first assertion of Lemma A(2) of the appendix,
$L_j^*\setminus A_j^*$ is an open $n$-disk, and hence the same 
holds true for $L_j\setminus A_j$.
Put $\Omega_j=B_{R_i+\epsilon}\cap B_{r_j}(p_{i,j},X)$ and note that
it is a convex set in $X$ containing $p_{i,j}$ in its interior.
It follows that the natural conical projection $c:L_j\to\partial\Omega_j$
(towards $p_{i,j}$) is a homeomorphism.
It is also not hard to realize that $C_j=c(A_j)$.

Now, choose $\rho\in(0,\pi)$ so close to $\pi$ that, putting 
$N_j^*:=L_j^*\setminus B^\bullet_\rho(w_j^*,L_j^*)$,
and denoting by $N_j$ the corresponding subset in $L_j$,
this subset $N_j$ is disjoint with the ball $B_{R_i+\epsilon}$,
and consequently its image $N'_j:=c(N_j)$ falls in
this part of $\partial\Omega_j$ which is
contained in $S_{R_i+\epsilon}$. By part (2) of Lemma A of the appendix, 
there is a collared $n$-disk $B^*$ in $L_j^*$ such that
$B_\rho(w_j^*,L_j^*)\subset B^*\subset B^\bullet_\pi(w_j^*,L_j^*)$.
Obviously, the corresponding $n$-disk $B$ in $L_j$ is also collared,
so we get that $L_j\setminus\hbox{int}(B)$ is homeomorphic  to
$M\setminus\hbox{int}(\Delta)$. We also have $L_j\setminus B\subset N_j$,
so that the image 
$P_j:=c(L_j\setminus B)$ is a subset of this part of $\partial\Omega_j$
which is contained in $S_{R_i+\epsilon}$.
Moreover, since $L_j\setminus A_j$ is an open $n$-disk,
it follows from the Annulus Theorem
that $[L_j\setminus\hbox{int}(B)]\setminus A_j$
is homeomorphic to the product $S^{n-1}\times[0,1)$.
Since the conical projection $c$ is a homeomorphism,
we also get that $P_j\cong M\setminus\hbox{int}(\Delta)$
and $P_j\setminus C_j\cong S^{n-1}\times[0,1)$.
Put $D'_j:=\psi_i'(P_j)$. By Lemma 4.3, 
$D_j'$ is homeomorphic to the quotient space $P_j/C_j$,
and by what was said above, this quotient is homeomorphic
to the closed $n$-disk which contains $p_{i,j}$ in its interior.
Clearly, we have 
$(\psi_i')^{-1}(D_j')=P_j\cong M\setminus\hbox{int}(\Delta)$.
By the Annulus Theorem, the same is true for any collared $n$-disk
$D$ contained in the interior of $D_j'$ and containing $p_{i,j}$
in its interior, hence the lemma.

\medskip\noindent
{\bf Lemma 4.5.}
{\it Let $g_i:S_{R_{i+1}}\to S_{R_i}$ be the geodesic projections, and 
$\Lambda_{R_i}=\{ p_{i,1},\dots,p_{i,k_i} \}$ be
the sets of singular points in the spheres $S_{R_i}$, as
in the proof of Main Theorem. Let $Z$ be an 
arbitrary finite subset of $S_{R_i}$
disjoint with $\Lambda_{R_i}$, and for each $z\in Z$
let $y_z\in S_{R_{i+1}}$ be an arbitrary point in the preimage 
$(g_i)^{-1}(z)$. Then $g_i$ can be approximated arbitrarily close
by a map $g_i^Z:S_{R_{i+1}}\to S_{R_i}$ satisfying conditions (2) and (3)
of Proposition 3.3 (with $g_i^Z$ substituted for $g_i'$) and such that}
$$
g_i^Z(y_z)=z \hbox{ \it for all } z\in Z.  \leqno{(C)}
$$

\medskip\noindent
{\bf Proof:} Choose positive $\epsilon$ as small as required
in Lemma 4.4, and consider the corresponding geodesic projections
$\psi_i:S_{R_i+\epsilon}\to S_{R_i}$ and 
$\phi_i:S_{R_{i+1}}\to S_{R_i+\epsilon}$ towards $x_0$. 
We then clearly have $g_i=\psi_i\phi_i$.
Consider a very close approximation $\psi_i'$ of $\psi_i$
satisfying conditions (1) and (2) of Lemma 4.3 and condition
($\star$) of Lemma 4.4. Moreover, since
by Lemma 3.2 the map $\phi_i$ is cell-like, 
consider its very close approximation by a homeomorphism $\phi_i'$,
as guaranteed by Approximation Theorem (recalled in the appendix).
The composition map $\psi_i'\phi_i'$ satisfies conditions (2) and (3)
of Proposition 3.3, but not necessarily the condition (C) of the assertion.

Note that, by the choices of $\psi_i'$ and $\phi_i'$, for each $z\in Z$
the point $\psi_i'\phi_i'(y_z)$ is very close to 
$\psi_i\phi_i(y_z)=g_i(y_z)=z$.
Since $S_{R_i}$ is a manifold,
we can choose a correcting homeomorphism 
$\omega:S_{R_i}\to S_{R_i}$, very close to the identity,
fixing all points of $\Lambda_{R_i}$,
and such that $\omega(\psi_i'\phi_i'(y_z))=z$ for all $z\in Z$.
The lemma follows by taking $g_i^Z:=\omega\psi_i'\phi_i'$,
since clearly such a composition can be chosen to approximate
the map $g_i$ arbitrarily close.

\medskip
We now pass to the final phase of the proof of Proposition 3.3.
Observe that Lemma 4.5 guarantees existence of approximations
of the maps $g_i$ satisfying conditions (2) and (3) of the proposition.
It remains to justify that appropriately chosen such approximations
satisfy also the last condition (4).

Chose any sequence $i_m$ of natural numbers such that $i_m\le m$
and each $i\ge1$ appears in this sequence infinitely often.
Choose also a sequence $\epsilon_m$ of positive reals converging to 0.
Proceed inductively with respect to $m$, as follows.
For $m=1$ and for each $M\in{\cal M}$ choose a finite subset
$Y_{M,1}\subset\Lambda_M\cap[X\setminus B_{R_1}]$ such that its image 
$Z_{M,1}=G_1(Y_{M,1})$ is an $\epsilon_1$-net in $S_{R_1}$
(i.e. any point of $S_{R_1}$ lies at distance at most $\epsilon_1$
from $Z_{M,1}$).
Put $Y_1=\cup_{M\in{\cal M}}Y_{M,1}$ and $Z_1=\cup_{M\in{\cal M}}Z_{M,1}$,
and assume, without loss of generality, that $Z_1$ is disjoint with
$\Lambda_{R_1}$. 
(The possibility of such a choice of the above sets follows from Claim 4.1.)
Put $\nu_1$ to be the largest $i$ such that
the intersection $Y_1\cap S_{R_i}$ is nonempty.
For $1\le i\le\nu_1-1$ put $Y_{1,i}:=Y_1\cap[X\setminus B_{R_i}]$,
$Z_{1,i}:=G_i(Y_{1,i})$, and for each $z\in Z_{1,i}$ put 
$y_z=G_{i+1}(y)$, where $y\in Y_{1,i}$ is this point for which 
$z=G_i(y)$.
Again for $1\le i\le\nu_1-1$ choose successively
approximations $g_i'=g_i^{Z_{1,i}}$ as in Lemma 4.5, 
so close to $g_i$ that the requirements
of Brown's Lemma are fulfiled.
Note that, apart from satisfying conditions (2) and (3)
of Proposition 3.3, these maps $g_i':1\le i\le\nu_1-1$ have the following
property: for each $M\in{\cal M}$ the images in $S_{R_1}$
of the points from $\cup_{i=1}^{\nu_1}(\Lambda_M\cap S_{R_i})$,
through the appropriate compositions of the maps $g_i'$,
form an $\epsilon_1$-net in $S_{R_1}$.

Now, assume that $m\ge2$, $\nu_{m-1}\ge m$, and that we have already defined
the approximations $g_i'$ for all $i\le\nu_{m-1}-1$.
For any $M\in{\cal M}$ choose a finite subset
$Y_{M,m}\subset\Lambda_M\cap[X\setminus B_{R_{\nu_{m-1}}}]$ such that its image 
$Z_{M,m}=G_{\nu_{m-1}}(Y_{M,m})$, further projected by the composition
$g_{i_m}'\circ g_{i_m+1}'\circ\dots\circ g_{\nu_{m-1}-1}'$,
is an $\epsilon_m$-net in $S_{R_{i_m}}$.
Put $Y_m=\cup_{M\in{\cal M}}Y_{M,m}$ and $Z_m=\cup_{M\in{\cal M}}Z_{M,m}$,
and assume, without loss of generality, that the image of $Z_m$ 
in $S_{R_{i_m}}$, through the above mentioned composition of the maps $g_i'$,
omits the points
$p_{i_m,1},\dots,p_{i_m,k_{i_m}}$. 
(The possibility of such a choice of the above sets follows from Claim 4.1,
and from the fact that the image of a dense subset through a surjective map
is dense.)
Put $\nu_m$ to be the largest $i$ such that
the intersection $Y_m\cap S_{R_i}$ is nonempty.
For $\nu_{m-1}\le i\le\nu_m-1$ put $Y_{m,i}:=Y_m\cap[X\setminus B_{R_i}]$,
$Z_{m,i}:=G_i(Y_{m,i})$,
and for each $z\in Z_{m,i}$ put $y_z=G_{i+1}(y)$ for this $y\in Y_{m,i}$
for which $G_i(y)=z$.
Then for $\nu_{m-1}\le i\le\nu_m-1$ choose successively
approximations $g_i'=g_i^{Z_{m,i}}$ as in Lemma 4.5, 
close enough to $g_i$ to fulfil the requirements
of Brown's Lemma.
Note that, apart from satisfying conditions (2) and (3)
of Proposition 3.3, the maps $g_i':1\le i\le\nu_m-1$ have the following
property: for each $M\in{\cal M}$ and for each $j\in\{ i_1,\dots,i_m \}$, 
the images in $S_{R_j}$
of the points from $\cup_{i=j+1}^{\nu_m}(\Lambda_M\cap S_{R_i})$, 
through the appropriate compositions of the maps $g_i'$,
form an $\epsilon_k$-net in $S_{R_j}$, where $k\le m$ is the largest
number such that $i_k=j$.

A direct verification shows that the whole sequence of maps $g_i':i\ge1$
described above meets all the requirements of Proposition 3.3,
which finishes the proof.

\bigskip\noindent
{\bf 5. Applications of Main Theorem.}

\medskip
In this section we describe vast classes of examples to which Main Theorem
(as presented in Section 2) applies. Among others, 
we provide proofs of Theorems 1 and 2
from the introduction.

\bigskip\noindent
{\it 
%5.1 
Hyperbolizations of $\cal M$-pseudomanifolds
and the proof of Theorem 1.}

\medskip
Given a finite family $\cal M$ of PL manifolds as in Section 2 
(all of the same dimension $n$),
{\it a closed $\cal M$-pseudomanifold} 
is a compact connected polyhedral cell complex $P$ such that
\item{(1)} for each point $x\in P$ the link $\hbox{Lk}(x,P)$
is either a PL $n$-sphere, or is PL-homeomorphic to some $M\in{\cal M}$;
furthermore, if all $M\in{\cal M}$ are oriented, we assume that $P$
is also oriented, and that each homeomorphism $\hbox{Lk}(x,P)\to M$
as above respects the orientations;
\item{(2)} for each $M\in{\cal M}$ the set 
$
\Lambda^P_M:=
\{ x\in P:\hbox{Lk}(x,P)\cong M \}$ is nonempty  
(if all $M\in{\cal M}$ are oriented, the symbol $\cong$ above
denotes relation of being PL homeomorphic as oriented manifolds).

\medskip
Recall that {\it hyperbolization}, as described e.g. in [DJ],
is a procedure that turns an arbitrary compact simplicial complex 
into a compact
nonpositively curved piecewise euclidean polyhedral complex
with the same local PL topology.
In particular, if we apply this procedure to a simplicial closed
$\cal M$-pseudomanifold $P$, we get a nonpositively curved 
piecewise euclidean complex $P_h$ which is also a closed
$\cal M$-pseudomanifold. (If all manifolds in $\cal M$ are oriented,
we consider the procedure which respects all requirements 
concerning orientations, e.g. the procedure described in
Subsection (4c) in [DJ], for which each hyperbolized simplex is an oriented
manifold, and the associated hyperbolization map is of degree 1.)
Observe that the universal cover
$\widetilde P_h$ is then an $({\cal M},0)$-pseudomanifold.

There is also an analogous procedure, called {\it strict hyperbolization}
and described in [CD], which turns any simplicial closed 
$\cal M$-pseudomanifold $P$ into a piecewise hyperbolic
locally $CAT(-1)$ closed $\cal M$-pseudomanifold $P_h^s$.
The universal cover
$\widetilde P_h^s$ is then an $({\cal M},-1)$-pseudomanifold.

\medskip
Theorem 1 of the introduction is an easy consequence of the following.

\medskip\noindent
{\bf 5.1 Theorem.}
{\it Let ${\cal M}=\{ M_1,\dots,M_k \}$ be a finite family of 
closed connected PL manifolds of the same dimension $n$,
either all oriented, or at least one of which is non-orientable.
Suppose that for some positive integers $m_1,\dots,m_k$
the disjoint union $\sqcup_{i=1}^k m_iM_i$ bounds a compact
$(n+1)$-dimensional PL manifold $W$ (in oriented sense,
if all manifolds from $\cal M$ are oriented). Then there is a hyperbolic
group $G$ such that its Gromov boundary $\partial G$ is homeomorphic
to the tree of manifolds ${\cal X}({\cal M})$.}

\medskip\noindent
{\bf Proof:}
Consider any PL triangulation of any manifold $W$ as in the assumptions. For each boundary component $M$ of $W$ consider the simplicial cone
over $M$ and glue it to $W$ via the identity of $M$. This gives a simplicial
closed $\cal M$-pseudomanifold which we denote by $P$. 
Consider its strict hyperbolization $P_h^s$, and its universal cover
$\widetilde P_h^s$, which is a $({\cal M},-1)$-pseudomanifold.
The group $G=\pi_1(P_h^s)$ is then a word hyperbolic group,
and its Gromov boundary $\partial G$ coincides with the visual boundary
$\partial \widetilde P_h^s$. Since by Main Theorem the latter boundary
is homeomorphic to the tree of manifolds ${\cal X}({\cal M})$,
this completes the proof.

\medskip\noindent
{\bf Proof of Theorem 1:}
Part (1) follows from Theorem 5.1 by taking ${\cal M}=\{M\}$,
while part (2) follows by taking ${\cal M}=\{N\}$ and $W=N\times[0,1]$.

\medskip\noindent
{\bf 5.2 Remark.}
Note that Theorem 1 provides new examples of Gromov boundaries
already in dimension 3. More precisely, since each closed connected
PL 3-manifold $M$ bounds a compact PL 4-manifold, each tree of
3-manifolds $M$, ${\cal X}(M)$, is homeomorphic to the Gromov
boundary of some hyperbolic group.
On the other hand, in the case of orientable 3-manifolds
the arguments of [PS] (after appropriate correction of the main
result of [F] used in [PS]) justify this statement only for
manifolds of form $M\#\overline{M}$, where $M$ and $\overline{M}$
are the oppositely oriented copies of any orientable 3-manifold.

\medskip\noindent
{\bf 5.3 Question.}
Theorem 5.1 leaves open the following question:
{\it is there a hyperbolic group $G$ whose Gromov boundary $\partial G$
is homeomorphic to the tree of complex projective planes
${\cal X}(CP^2)$?} Recall that neither $CP^2$ nor any positive number
of its copies bounds a compact oriented 5-manifold, so $G$ cannot be
obtained by referring to Theorem 5.1. On the other hand, this theorem
easily implies that the space ${\cal X}(\{ CP^2,\overline{CP^2} \})$
is homeomorphic to the Gromov boundary of some hyperbolic group.
Since one can use properties of \v Cech cohomology rings to show that
the spaces ${\cal X}(CP^2)$ and ${\cal X}(\{ CP^2,\overline{CP^2} \})$
are not homeomorphic, this also does not help to answer the question.
My guess is that the answer is negative.
Of course, a similar question can be asked for other than $CP^2$
manifolds which represent the elements of infinite order
in the corresponding oriented cobordism additive semi-group.

\bigskip\noindent
{\it Coxeter groups with manifold nerves and the proof of Theorem 2.}

\medskip
Our main reference concerning Coxeter groups is the book [Da] by Mike
Davis. We start with explaining the terms appearing in the statement
of Theorem 2, following the terminology and notation from [Da].

Given a finite set $S$, a {\it Coxeter matrix} on $S$ is a matrix
${\bf m}=(m_{st})_{s,t\in S}$ such that
\item{(1)}  for each $s\in S$ we have $m_{ss}=1$;
\item{(2)} for all $s,t\in S$, $s\ne t$, we have that $m_{st}$ is
an integer $\ge2$ or $\infty$, and $m_{st}=m_{ts}$.

\noindent
Associated to a Coxeter matrix $\bf m$ on $S$, there is a group $W$
given by the presentation
$$\langle  S\,|\,\{ (st)^{m_{st}}:s,t\in S \}  \rangle,$$
where the symbol $(ab)^\infty$ denotes absence of any relation
of the form $(ab)^k$ in the set of relations. $W$ is called
the {\it Coxeter group} associated to $\bf m$, and it is known
that the set $S$ canonically injects in $W$. The pair $(W,S)$ is 
called the {\it Coxeter system} associated to $\bf m$.

For any subset $T\subset S$ the {\it special subgroup} $W_T<W$
is the subgroup generated by $T$. It is known that $(W_T,T)$
can be canonically identified with the Coxeter system associated
to the restricted matrix ${\bf m}_T$, i.e. the matrix $(m_{st})_{s,t\in T}$.

The {\it nerve} of a Coxeter system $(W,S)$ is the simplicial complex
$L=L(W,S)$ whose vertex set coincides with $S$, and such that
$T\subset S$ spans a simplex of $L$ iff the special subgroup $W_T$
is finite.

As it is described in Chapter 7 of [Da], to any Coxeter system $(W,S)$
there is associated a polyhedral cell complex $\Sigma=\Sigma(W,S)$,
called the {\it Davis-Moussong complex} of $(W,S)$, which satisfies
the following properties:
\itemitem{($\Sigma1$)} each vertex link of $\Sigma$ is a simplicial
complex isomorphic with the nerve $L$;
\itemitem{($\Sigma2$)} $\Sigma$ carries a natural piecewise euclidean
metric with respect to which it is a $CAT(0)$ space,
see Theorem 12.3.3 on p. 235 in [Da];
\itemitem{($\Sigma3$)} the group $W$ acts on $\Sigma$ by isometries,
properly discontinuously and cocompactly, so that the generators
from $S$ correspond to certain geometric reflections in $\Sigma$,
and so that the action is simply transitive on the vertex set of $\Sigma$.

\medskip\noindent
{\bf Proof of Theorem 2:}
Suppose that the nerve $L(W,S)$ is a PL triangulation of a closed
connected manifold $M$. If $M$ is orientable, it follows from conditions 
($\Sigma1$)-($\Sigma3$) above that $\Sigma(W,S)$ is an
$({\cal M},0)$-pseudomanifold with ${\cal M}=\{ M,\overline M \}$,
where $M$ and $\overline M$ are the two oppositely oriented
copies of $M$. Similarly, if $M$ is non-orientable, $\Sigma(W,S)$ is an
$({\cal M},0)$-pseudomanifold with ${\cal M}=\{ M \}$.
Thus, applying Main Theorem to the Davis-Moussong complex
$\Sigma(W,S)$, we immediately get Theorem 2.

\medskip\noindent
{\bf 5.4 Remark.}
It is known that, when a Coxeter group $W$ is word hyperbolic,
its Gromov boundary $\partial W$ coincides with the visual
boundary $\partial\Sigma(W,S)$ (see e.g. Remark I.8.5 on p. 527 in [Da]).
In such a case, if the nerve $L(W,S)$
is a PL triangulation of a closed connected manifold $M$,
the Gromov boundary $\partial W$ is homeomorphic to the
space ${\cal X}({\cal M})$ as in the statement of Theorem 2.

\bigskip\noindent
{\bf 6. Riemannian $({\cal M},0)$-pseudomanifolds with log-injective
singularities.}

\medskip
In this section we explain how to adapt our proof of Main Theorem
to a class of smooth $CAT(0)$ pseudomanifolds with Riemannian
metrics on their regular part. Using this, we deduce Theorem 3.
We start with describing the relevant class of pseudomanifolds.
We refer the reader to Subsection 3.2 of [FM] for the definition of
the {\it space of directions} $\Sigma_p(X)$ at a point $p$ of a $CAT(0)$
space $X$ (or to Definition II.3.18 on p. 190 of [BH], where
the same object is denoted $S_p(X)$).

\medskip\noindent
{\bf Definition 6.1.}
Given a finite collection $\cal M$ of closed connected manifolds
of the same dimension $n$, a $CAT(0)$ complete geodesic
metric space $(X,d)$ is a {\it Riemannian $({\cal M},0)$-pseudomanifold
with log-injective singularities} if the following conditions are satisfied:

\smallskip
\item{(1)} there is a subset $\Lambda\subset X$, called {\it the singular
set} of $X$, which is discrete and the complement $X\setminus\Lambda$
is a smooth manifold;

\item{(2)} the set $\Lambda$ is partitioned into subsets 
$\Lambda_M:M\in{\cal M}$ such that for each $M\in{\cal M}$
and any $p\in\Lambda_M$ the space of directions $\Sigma_p(X)$ 
is homeomorphic to $M$;

\item{(3)} each of the subsets $\Lambda_M$ is a net in $X$;

\item{(4)} the {\it regular part} $X\setminus\Lambda$ is equipped with
a Riemannian metric $g$ of nonpositive sectional curvature, such that
$d$ restricted to $X\setminus\Lambda$ coincides with the path metric
induced by $g$, and the completion of $(X\setminus\Lambda,d)$
coincides with $(X,d)$;

\item{(5)} each $p\in\Lambda$ has a neighborhood $U$ 
such that the logarithmic map 
$\log_p:U\setminus\{p\}\to\Sigma_p(X)\times R_+$
(which to each $x\in U\setminus\{ p \}$ associates the pair
$(a,r)$ such that $a\in\Sigma_p(X)$
is the direction of the geodesic $[p,x]$ and $r=d(p,x)$) is injective;

\item{(6)} the space $\Sigma_p(X)$ at any singular point $p$ has the
property that every ball of radius $r\in(0,\pi)$ in it
(with respect to the angle metric) is a collared $n$-disk in 
$\Sigma_p(X)$.

\medskip\noindent
{\bf Theorem 6.2.}
{\it Let $X$ be a Riemannian $({\cal M},0)$-pseudomanifold
with log-injective singularities. Then the visual boundary $\partial X$
is homeomorphic to the tree of manifolds ${\cal X}({\cal M})$.}

\medskip
Theorem 6.2 follows by the same arguments as in the proof of
Main Theorem, slightly adapted and simplified according to the
following features:

\smallskip
\item{(1)} whenever in the proof of Main Theorem we use
cone neighbourhoods of singular points,
we need to use instead small balls which are log-injective neighbourhoods
of singular points; existence of the latter balls is justified by condition (5)
in Definition 6.1;

\item{(2)} references to Lemma A of the appendix appearing 
in the proof of Main Theorem,
when applied to links at singular points, need to be replaced by
references to  condition (6) of Definition 6.1; 
moreover, references to properties of links at non-singular points
also need to be replaced by references to the properties of
the corresponding spaces of directions, which are just the
standard round $n$-spheres of constant curvature 1;

\item{(3)} since geodesics in 
Riemannian $({\cal M},0)$-pseudomanifolds
do not bifurcate outside the singular set, while proving Theorem 6.2
we never need to approximate 
(by referring to Approximation Theorem for cell-like maps)
various geodesic projections between
the spheres (or their restrictions) by homeomorphisms,
since in this setting the corresponding maps are automatically homeomorphisms.

\medskip\noindent
{\bf Proof of Theorem 3:}
In the paper [FM] by K. Fujiwara and J. Manning
it is explained how to put a Riemannian metric  on a regular part
of any space  appearing 
in the statement of Theorem 3, so that its lift to the universal
cover of this space (and the induced path metric and its completion)
satisfies all the requirements of Definition 6.1.
In fact, the metrics constructed in [FM] are even {\it CAT}$(-1)$.
In view of this, Theorem 3 follows fairly directly from Theorem 6.2.

\medskip
\noindent
{\it A related class of examples and an open question.}

\smallskip
In [C] Coulon studies a class of examples related to those
appearing in Theorem 3. Namely, he considers a negatively
curved closed Riemannian manifold $M$ and its closed totally
geodesic submanifold $N$, and the space 
$Y=M\cup_N\hbox{Cone}(N)$ obtained by attaching to $M$
a cone of base $N$.
He shows that, under certain conditions,  the resulting space
is aspherical, and its fundamental group $\Gamma=\pi_1Y$ is hyperbolic.

Note that the same group $\Gamma$ is the fundamental group
of the pseudomanifold $X$ obtained from $M$ by simply
collapsing $N$ to a point. Equivalently, we can delete
from $M$ some open tubular neighbourhod $V$ of $N$,
and collapse the boundary $\partial V=\partial(M\setminus V)$
to a point, which we denote $p$. Clearly, $X$ is then 
homotopy equivalent to $Y$. Moreover, it is
a pseudomanifold with natural PL structure, and $p$ is its
only singular point. The link $\hbox{Lk}(p,X)$
is PL homeomorphic to $\partial V$, and thus $X$ is a
$\{ \partial V \}$-pseudomanifold (as defined at the beginning 
of Section 5).

The arguments used by Coulon in [C] do not guarantee
existence of a negatively or nonpositively curved metric on $X$.
Thus, the following question naturally appears.

\medskip\noindent
{\bf Question 6.3.} Is Gromov boundary $\partial\Gamma$
of a hyperbolic group $\Gamma$ as above, resulting from the
construction of Coulon [C], homeomorphic to the tree
of manifolds ${\cal X}(\partial V)$?

\bigskip
\noindent
{\bf Appendix.} 

\medskip
In this appendix we present (slightly changed and extended variants of) some results
from Section 3 of the paper [DJ] by M. Davis and T. Januszkiewicz
(see Lemmas A--C below). 
These results are necessary ingredients in the arguments in Sections 3 and 4 
of the present paper.
We provide also proofs of these results. 
The reasons for doing this are
as follows:
\item{(a)} some of these results are not explicitely stated in [DJ], 
though they appear implicitely in the proof of Lemma 3b.1 and Theorem 3b.2
on pages 371--372 of that paper; 
\item{(b)} the statements we give are more general than those in [DJ],
although they follow by essentially the same arguments;
for our applications in Sections 3 and 4
we need the extended versions of the results;
\item{(c)} there are two gaps in the above mentioned proof of Lemma 3b.1
and Theorem 3b.2 in [DJ], as it was pointed to me by the referee
of the present paper; the first gap was found by Hanspeter Fischer
and rectified in a 1997 letter of Mike Davis to Fischer,
but the new argument has been never published; 
the idea how to rectify the second gap was communicated to me
by the referee of the present paper;
the essential part
of this appendix is devoted to presentation of these two new arguments; 
see Remark 2 below for more details; 
%%concerning the mistake,
%%and for a precise indication of the part of the appendix corresponding
%%to the new argument provided by Mike Davis;
\item{(d)} the proofs of (the analogons of) Lemmas A--C
given in [DJ] form a scheme of induction (with respect to the
dimension $n$) involving all of these results at once; 
even though the first above mentioned mistake concerns the proof
of Lemma C, and second one the proof of Lemma A,
it is hard to present the correction without mentioning the other results,
and providing the whole of their common proof.  

\medskip
We thank the referee of the present paper for his suggestions concerning
the scope and the shape of this appendix.

\bigskip
\noindent
{\it The results.}
\medskip
We use the terminology as in [BH]. In our notation,
$B_r(x,X)$ is the {\it closed} metric ball, and $B_r^\bullet(x,X)$
the {\it open} one.
The notions of a cellular set  
and a cell-like set and map (appearing in the statements of Lemmas A and C)
are recalled in the next part of this appendix.

\medskip
\noindent
{\bf Lemma A.}
{\it Let $L$ be a $CAT(1)$ piecewise spherical closed PL manifold of dimension $n$,
and let $v\in L$ be any point.} 
\item{(1)} 
{\it For any $r\in(0,\pi/2)$ the closed ball $B_r(v,L)$
is homeomorphic to the $n$-disk $B^n$, with the sphere $S_r(v,L)$
corresponding to the boundary $\partial B^n$, and 
it is collared in $L$.}
\item{(2)}
{\it For any $r\in(0,\pi]$ the open ball $B^\bullet_r(v,L)$ is homeomorphic to the
open $n$-disk. As a consequence, for any $r\in(0,\pi)$ and any $\varepsilon\in(0,\pi-r)$
there is a collared $n$-disk $D$ in $L$ such that 
$B_r(v,L)\subset D\subset B^\bullet_{r+\varepsilon}(v,L)\subset B^\bullet_\pi(v,L)$.}
\item{(3)}
{\it If moreover $L$ is the standard PL $n$-sphere, 
then for any $r\in(0,\pi)$
the complement
$L\setminus B^\bullet_r(v,L)$ is a cellular subset of $L$.
In particular,  $L\setminus B^\bullet_\pi(v,L)$ is cellular in $L$.}

\medskip
\noindent
{\bf Lemma B.}
{\it Let $X$ be a locally compact  
$CAT(\kappa)$ 
$M_\kappa$-polyhedral complex with $\hbox{Shapes}(X)$ finite.
Fix any $x_0\in X$ and any $r>0$ ($r<\pi/2\sqrt\kappa$ if $\kappa>0$).
Suppose that for each $x\in X$ such that $d(x,x_0)=r$
the metric link
$\hbox{\rm Lk}(x,X)$ (viewed as a piecewise spherical complex) is
a closed PL manifold of dimension ${n-1}$. Then the sphere $S_r(x_0,X)$ is a closed $(n-1)$-manifold,
and a neighbourhood $U$ of $S_r(x_0,X)$ in $B_r(x_0,X)$
is an $n$-manifold with boundary $\partial U=S_r(x_0,X)$.}

\medskip
\noindent
{\bf Lemma C.}
{\it Let $X$ be a locally compact 
$CAT(\kappa)$ 
$M_\kappa$-polyhedral complex with $\hbox{Shapes}(X)$ finite.
Let $x_0\in X$, and for any
$t>r>0$ such that $t<\pi/\sqrt{\kappa}$ if $\kappa>0$ 
denote by $c_{t,r}:S_t(x_0,X)\to S_r(x_0,X)$
the geodesic projection towards $x_0$ between the metric spheres.
Suppose that for each $x\in X$ such that $r\le d(x,x_0)<t$ the metric link
$\hbox{\rm Lk}(x,X)$ (viewed as a piecewise spherical complex) is PL homeomorphic
to the standard PL sphere $S^{n-1}$. Then $c_{t,r}$ is a cell-like map.}

\medskip
Before passing to the proof of the above lemmas we need some 
preparations.

\bigskip
\noindent
{\it Preparations concerning cellular sets and cell-like sets and maps.}
\medskip

A nonempty compact subset $C$ of a metric $n$-manifold  $(M,d)$ is {\it cellular}
if $C$ has arbitrarily close open neighbourhoods which are homeomorphic
to the open $n$-disk. 
More precisely, given any $\epsilon>0$
there is a neighbourhood of $C$ of the relevant form mentioned above
which is contained in the $\epsilon$-neighbourhood
$N_\epsilon(C)=\{ x\in M:d(x,C)<\epsilon \}$.

A nonempty compact metric space $C$ is {\it cell-like} if it can be
embedded into the Hilbert cube $Q$ so that for any neighbourhood $U$
of $C$ in $Q$ the space $C$ is null-homotopic in $U$.
It is known (see e.g. [Ed], remark on p. 114)  that 
a finite dimensional compact metric space is cell-like if it can be embedded
as a cellular subset in some manifold. Obviously,
each cellular subset of a manifold is cell-like, but the converse
is not true (i.e. there are embeddings of cell-like spaces
in manifolds with non-cellular images). A map between metric spaces is {\it cell-like}
if it is a proper surjection and each point preimage is cell-like.

\medskip
The following result is due to Siebenmann [Si] for $n\ge5$,
Quinn [Qu] for $n=4$, Armentrout [Ar] for $n=3$ and R. L. Moore
[Mo] for $n\le2$. The case of manifolds with nonempty boundary
is carefully addressed in [Si], and the cases involving dimension 3 hold
because the Poincare conjecture is true.

\medskip\noindent
{\bf Approximation Theorem.}
{\it Each cell-like map between manifolds with boundary 
is a near-homeomor\-phism,
i.e. it can be approximated by homeomorphisms. More precisely,
if $f:M\to N$ is cell-like and $d$ is any metric in $N$ (compatible
with the topology) then for any continuous positive function
$\delta:M\to R_+$ there is a homeomorphism $h:M\to N$
such that for all $x\in M$ we have $d(f(x),h(x))\le\delta(x)$.}

%\medskip\noindent
%{\bf Remark 1.}
%Observe that for any near-homeomorphism between the manifolds,
%any point preimage is easily seen to be a cellular subset. Thus,
%by Approximation Theorem, a cell-like map between manifolds is 
%automatically cellular.

\medskip
The proof of the next result was indicated to me by the referee
of the present paper. 

\medskip\noindent
{\bf Lemma 1.}
{\it If $f:A\to B$ is a cell-like surjection between finite dimensional
metric compacta, then $A$ is cell-like if and only if $B$ is cell-like.}

\medskip\noindent
{\bf Proof:} 
Recall that ENR is a locally compact finite dimensional ANR.
Embed $A$ in a compact ENR $X$
(e.g. a sphere of dimension $2\dim A+1$). Let $Y$ be the quotient space
obtained from $X$ by shrinking each set $f^{-1}(b):b\in B$ (viewed 
as a subset of $X$) to a point. Since $Y$ is easily seen to be a Hausdorff
space, the quotient map  $F:X\to Y$ is closed, and hence also perfect. 
Since metrizability is preserved by perfect maps (see e.g. [En],
Theorem 4.4.15), $Y$ is metrisable.

By definition of $Y$, $F$ is a cell-like map.
Since $\dim Y\le\dim B+\dim(X\setminus A)+1$,
we get that $\dim Y$ is finite.
It follows then from Corollary 3.3 in [La] that $Y$ is an ENR.
Since $F^{-1}(B)=A$, Theorem 1.4 of [La] implies that
$A$ is cell-like if and only if $B$ is cell-like, as required.

\bigskip
\noindent
{\it The scheme of the inductive proof of Lemmas A--C.}
\medskip

Observe that for $n=1$ the statements of Lemmas A--C
are obviously true. Denote by $(A_n)$, $(B_n)$
and $(C_n)$ the statements of Lemmas A--C in dimension $n$,
respectively. We shall prove Lemmas A--C according to the 
scheme
$$
(A_{n-1})\Rightarrow(B_n)\Rightarrow(C_n)\Rightarrow(A_n).
$$

\medskip\noindent
{\bf Remark 2.} 
The first gap in the paper [DJ], mentioned in comment (c)
at the beginning of this appendix, concerns the argument for showing
the analogue of our
implication $[(A_{n-1})\hbox{ and }(B_n)]\Rightarrow(C_n)$.
The gap appears in the paragraph beginning 
in the middle of page 372 of [DJ], and more precisely in the 
seventh sentence of this paragraph starting with "Assume that we 
have chosen $s$ close enough to $r$...". 
Examples can be constructed in which for every $s>r$ there are 
geodesic rays emanating from $x$ that have two or more bifurcation points
(i.e. points with nontrivial infinitesimal shadow) in $A_{r,s}$.
The new argument provided by the authors of [DJ], 
which rectifies this gap,
is presented below, in the part concerning the proof of the implication
$[(A_{n-1})\hbox{ and }(B_n)]\Rightarrow(C_n)$.

The second gap in [DJ] concerns the proof of Lemma (3b.1)
from that paper, and it appears in line -6 on p. 371.
In fact, the proof of the implication $(L_{n-1})\Rightarrow(T_n)$,
as presented in [DJ], does not apply to the implication 
$(L_{n-1})\Rightarrow(L_n)$, for balls with radius $r\in[\pi/2,\pi)$.
The reason is that balls with such radii need not be locally geodesically
strictly convex in the corresponding space $L$
(in the sense that the geodesic between any two sufficiently close
points on the boundary sphere of such a ball intersects this sphere
only at the endpoints). The latter property
is esssentailly exploited in the argument in [DJ], see e.g. the statement
in line  1 on p. 372 of that paper, and few next sentences.
In fact, it is not clear if Lemma (3b.1) from [DJ] holds true in its full
generality. Thus, in our exposition we only formulate (as Lemma A(1)) its
special case, for $r\in(0,\pi/2)$. We also prove some result slightly weaker
than Lemma (3b.1) of [DJ], namely Lemma A(2). We provide arguments
that replace those from [DJ], whenever the latter arguments refer
to the full statement of Lemma (3b.1).

\bigskip
\noindent
{\it Proof of $(A_{n-1})\Rightarrow(B_n)$.}
\medskip
The proof of the implication $(L_{n-1})\Rightarrow(T_n)$ in [DJ],
in a long paragraph beginning at the end of page 371 and ending
in the middle of the next page, works without essential changes
in a slightly more general setting of $(B_n)$. 
More precisely, for any point $y\in S_r(x_0,X)$,
denoting by $v\in\hbox{Lk}(y,X)$ the point induced by the geodesic segment
$[y,x_0]$,
the above mentioned argument from [DJ] shows that a small closed neighbourhood
of $y$ in $B_r(x_0,X)$ is homeomorphic to the 
space obtained from the product $B_{\pi/2}(v,\hbox{Lk}(y,X))\times[0,\varepsilon]$
by collapsing to a point its subset 
$$[B_{\pi/2}(v,\hbox{Lk}(y,X))\times\{0\}]\cup [S_{\pi/2}(v,\hbox{Lk}(y,X))
\times[0,\varepsilon)],
$$ 
where the latter point corresponds to $y$.
Since, by part (2) of $(A_{n-1})$, the open ball $B^\bullet_{\pi/2}(v,\hbox{Lk}(y,X))$
is homeomorphic to the open $(n-1)$-disk, it follows easily that 
the above quotient is homeomorphic to the closed $n$-disk $B^n$,
with the point corresponding to $y$ lying on its boundary $\partial B^n$. 
We omit further details.

%\medskip\noindent
%{\it Remark.} If $\kappa>0$, the above mentioned argument from [DJ]
%does not work for balls with radius $r\ge\pi/2\sqrt\kappa$, since such balls are not
%necessarily locally strictly convex (in the sense that the geodesic segment
%connecting any two sufficiently close points from the boundary sphere
%intersects this sphere only at the endpoints). 

\bigskip
\noindent
{\it Proof of $[(A_{n-1})\hbox{ and }(B_n)]\Rightarrow(C_n)$.}
\medskip

The argument below essentially coincides with the one communicated
by Mike Davis in his 1997 letter to Hanspeter Fischer.

\smallskip
Suppose $g:[a,b]\to X$ is a geodesic segment in a metric space $X$,
and let $a<t<b$. Recall that $g(t)$ is a {\it bifurcation point} of $g$
if there is a geodesic segment $g':[a,c]\to X$, with $t<c$, such that
$g(s)=g'(s)$ for all $a\le s\le t$ and $g(s)\ne g'(s)$ for all
$t<s\le\min(b,c)$. If $X$ is a $CAT(\kappa)$ $M_\kappa$-polyhedral
complex, and if $v$ is the point in the metric link $L=\hbox{\rm Lk}(g(t),X)$
corresponding to the incoming direction of $g$ (i.e. direction induced by
$g|_{[a,t]}$) then $g(t)$ is a bifurcation point of $g$ precisely when
the infinitesimal shadow of $g(t)$ with respect to $v$ (as defined in [DJ]
on page 369) is nontrivial, i.e. consists of more than one point.
This infinitesimal shadow, denoted $\hbox{\rm Shad}_{g(t)}(v)$, is known to
coincide with the set $\{ y\in L:d_L(y,v)\ge\pi \}$,
where $d_L$ is the canonical piecewise spherical metric in $L$.

\medskip\noindent
{\bf Fact 3.}
{\it If $g:[a,b]\to X$ is a geodesic segment in an $M_\kappa$-polyhedral
complex with $\hbox{\rm Shapes}(X)$ finite, then the set of bifurcation
points of $g$ is finite.}

\medskip\noindent
{\bf Proof:} Since each $M_\kappa$-polyhedral complex with finite shapes
can be subdivided into an isometric $M_\kappa$-simplicial complex
with finite shapes ([BH], Proposition 7.49, page 118), we can assume
$X$ is an $M_\kappa$-simplicial complex with finite shapes.
Since every geodesic segment in such a simplicial complex is 
the concatenation of a finite number of segments, each of which is 
contained in a single simplex ([BH], Corollary 7.29, page 110),
this also holds for $g$. If $\gamma$ is one of such subsegments for $g$,
we claim that no interior point $y$ of $\gamma$ is a bifurcation point of $g$.
Indeed,  since the interior of $\gamma$ is contained in the interior
of some simplex of $X$, if $v,w$ denote the points of the link
$\hbox{Lk}(y,X)$ corresponding to the incoming and outgoing direction
of $g$, then the whole metric link  $\hbox{Lk}(y,X)$ is isometric
to some spherical suspension in which $v$ and $w$ are the suspension points.
It follows that $w$ is the only point in $\hbox{Lk}(y,X)$ at distance $\pi$
from $v$, so the infinitesimal shadow of $y$ at $v$ is trivial.
Thus only the ends of the subsegments $\gamma$ of $g$ as above
can be the bifurcation points of $g$, hence the fact.

\medskip\noindent
{\bf Fact 4.}
{\it Let $X$ be a $CAT(\kappa)$ $M_\kappa$-polyhedral complex with
$\hbox{Shapes}(X)$ finite. Fix any $x_0\in X$, $\rho>0$ ($\rho<\pi/\sqrt\kappa$
if $\kappa>0$), $x\in S_\rho(x_0,X)$, and let $v$ be the point of the link
$\hbox{\rm Lk}(x,X)$ corresponding to the incoming direction
of the geodesic from $x_0$ to $x$. Then there is $\epsilon>0$ such that
if $\rho<s<\rho+\epsilon$ then the preimage $c_{s,\rho}^{-1}(x)$ is 
homeomorphic to the infinitesimal shadow $\hbox{\rm Shad}_x(v)$.}

\medskip\noindent
{\bf Proof:}
There is $\epsilon>0$ such that the open ball $B^\bullet_\epsilon(x,X)$
is naturally isometric to the open ball of radius $\epsilon$ centered at the cone 
point in the $\kappa$-cone $C_\kappa(\hbox{\rm Lk}(x,X))$
(see Theorem I.7.16 on page 103 in [BH]).
We may additionally assume that $\rho+\epsilon<\pi/\sqrt\kappa$ if
$\kappa>0$. It follows that for any $s$ with $\rho<s<\rho+\epsilon$
there is a canonical map $\Phi:\hbox{\rm Lk}(x,X)\to S_{s-\rho}(x,X)$,
which is a homeomorphism. Moreover, it is not hard to observe that $\Phi$
maps the infinitesimal shadow $\hbox{\rm Shad}_x(v)$ to the set
$c_{s,\rho}^{-1}(x)$, thus establishing a homeomorphism as required.

\medskip
We now pass to the essential part of the proof of the implication
$[(A_{n-1})\hbox{ and }(B_n)]\Rightarrow(C_n)$.  
Note that, due to the assumption concerning links,
the region $\{ z\in X:r< d(z,x_0)< t \}$ has the
geodesic extesion property (compare {BH, Proposition II.5.10, page 208}),
and thus $c_{t,r}$ is a surjection.
Moreover, since by our assumptions the spheres in $X$ are compact,
this map is also proper
Arguing by contradiction,
suppose that $c_{t,r}$ is not cell-like. Hence there is $x\in S_r(x_0,X)$ such that
$c_{t,r}^{-1}(x)$ is not a cell-like set. We will construct a geodesic
segment, starting at $x_0$ and passing through $x$, 
that contains an infinite sequence $(y_i)$ of bifurcation points,
thus contradicting Fact 3.

\medskip\noindent
{\bf Claim 1.}
{\it There is a point $y_1\in X$ such that} 
\item{(1)} {\it $x$ lies on the geodesic $[x_0,y_1]$ from $x_0$ 
to $y_1$, and}
\item{(2)} {\it every geodesic from $x_0$ to a point of $c_{t,r}^{-1}(x)$
passes through $y_1$ (i.e. contains $[x_0,y_1]$) and bifurcates at $y_1$.}

\medskip\noindent
{\bf Proof:}
Since the set $c_{t,r}^{-1}(x)$ is not cell-like, it contains more than 
one point. For any two distinct points of $c_{t,r}^{-1}(x)$ there is
a point $z$ where the two geodesics from these points to $x_0$
meet for the first time. Put $s=d(x_0,z)$, note that $r\le s<t$, and
consider the infimum $t_1$ of all such $s$, for all pairs of distinct points
in $c_{t,r}^{-1}(x)$. Then $r\le t_1<t$, and all geodesics from $x_0$
to a point of $c_{t,r}^{-1}(x)$ must coincide on the interval $[0,t_1]$.
Let $y_1\in S_{t_1}(x_0,X)$ be the common point that all these geodesics 
pass through. Then $y_1$ is a bifurcation point for all of these geodesics.
Indeed, if this is not true, the infinitesimal shadow at $y_1$ for all of
these geodesics must be trivial. The existence of a cone neighbourhood
of $y_1$ in $X$ yields then existence of $y'$ such that $y_1$ is an interior
point of the geodesic $[x_0,y']$ and all geodesics from $x_0$ to a point
of $c_{t,r}^{-1}(x)$ contain $[x_0,y']$. But this contradicts the choice
of $t_1$, thus completing the proof.

\medskip\noindent
{\bf Claim 2.}
{\it There is $r_1\in(t_1,t)$ and a point $x_1\in S_{r_1}(x_0,X)$ such that}
\item{(1)} {\it the points $x$ and $y_1$ lie on the geodesic $[x_0,x_1]$, and}
\item{(2)} {\it the set $c_{t,r_1}^{-1}(x_1)$ is not cell-like.}

\medskip\noindent
{\bf Proof:}
Consider the link $L=\hbox{\rm Lk}(y_1,X)$ and note that, by 
the assumptions in Lemma C, $L$ is a $CAT(1)$ piecewise spherical complex
(with $\kappa=1$) PL homeomorphic to the standard PL sphere $S^{n-1}$.
Let $v\in L$ be the point corresponding to the incoming direction
of the geodesic $[x_0,y_1]$. By part (3) of $(A_{n-1})$, the infinitesimal shadow
$\hbox{\rm Shad}_{y_1}(v)=L\setminus B_\pi^\bullet(v,L)$ is a cellular
subset of $L$, and hence it is cell-like. Moreover, it follows from Fact 4
that there is $r_1\in(t_1,t)$ such that the preimage $c_{r_1,t_1}^{-1}(y_1)$
(which coincides with the preimage $c_{r_1,r}^{-1}(x)$)
is homeomorphic to the infinitesimal shadow $\hbox{\rm Shad}_{y_1}(v)$,
and hence it is also a cell-like set. Note that for the geodesic projection
$c_{t,r_1}$ the preimage $c_{t,r_1}^{-1}(c_{r_1,r}^{-1}(x))$  coincides with the set $c_{t,r}^{-1}(x)$. 
Moreover, the restriction $f:c_{t,r}^{-1}(x)\to c_{r_1,r}^{-1}(x)$
of $c_{t,r_1}$ is a surjection, 
because due to the assumption concerning links,
we have the appropriate
geodesic extesion property.
Since the set $c_{t,r}^{-1}(x)$ is not cell-like,
it follows from Lemma 1 that there is $x_1\in c_{r_1,r}^{-1}(x)$ 
such that $f^{-1}(x_1)=c_{t,r_1}^{-1}(x_1)$ is not a cell-like set,
which completes the proof.

\medskip
Iterating Claims 1 and 2 we construct an infinite sequence of real numbers
$$
r\le t_1<r_1\le t_2<r_2\le t_3<r_3 \dots <t
$$
and sequences of points 
$$
y_i\in S_{t_i}(x_0,X) \hbox{\quad and \quad} x_i\in S_{r_i}(x_0,X)
\hbox{\quad for } i\ge1
$$
such that for each $i\ge1$ the geodesic $[x_0,x_i]$ bifurcates at $y_i$
and is contained in $[x_0,x_{i+1}]$.
By the latter property, there is a limit $x_*=\lim_{i\to\infty}x_i\in X$,
and the geodesic $[x_0,x_*]$ contains the union of the geodesics
$[x_0,x_i]$. In particular, all points $y_i$ belong to $[x_0,x_*]$,
and each of them is a bifurcation point for this geodesic.
Since this contradicts Fact 3, the assertion of $(C_n)$ follows.

\bigskip
\noindent
{\it Proof of $[(B_n) \hbox{ and }(C_{n})]\Rightarrow(A_n)$.}

\medskip
First, assuming $(B_n)$ and $(C_n)$, we establish two auxiliary facts.

\medskip
\noindent
{\bf Fact 5.}
{\it Under assumptions as in Lemma C, 
and assuming additionally that $t<\pi/2\sqrt\kappa$ if $\kappa>0$,
consider the annular region
$$
A_{r,t}=B_t(x_0,X)\setminus B^\bullet_r(x_0,X)=
\{ x\in X:r\le d(x,x_0\le t) \}.
$$
Then $A_{r,t}$ is homeomorphic to the product
$S_r(x_0,X)\times[r,t]$, with $S_r(x_0,X)\times\{r\}$ and 
$S_r(x_0,X)\times\{t\}$ corresponding to the spheres
$S_r(x_0,X)\subset A_{r,t}$ and $S_t(x_0,X)\subset A_{r,t}$,
respectively.}

\medskip\noindent
{\bf Proof:}
We will construct a homeomorphism
$A_{r,t}\to S_r(x_0,X)\times[r,t]$ which maps the spheres $S_t(x_0,X)$
and $S_r(x_0,X)$ on the sets $S_r(x_0,X)\times\{t\}$
and $S_r(x_0,X)\times\{r\}$, respectively.

Note first that $A_{r,t}\setminus S_r(x_0,X)$ is a manifold
with boundary $\partial(A_{r,t}\setminus S_r(x_0,X))=S_t(x_0,X)$. 
Indeed, 
due to the assuption on links,
each point in $A_{r,t}\setminus(S_r(x_0,X)\cup S_t(x_0,X))$
has a neighbourhood in $A_{r,t}\setminus S_r(x_0,X)$ 
homeomorphic to the open $n$-disk.
Furthermore, due to the last assertion in $(B_n)$, each point of $S_t(x_0,X)$
has a neigbourhood in $A_{r,t}$ homeomorphic to the open half-$n$-disk.

Consider now the map $\psi:A_{r,t}\setminus S_r(x_0,X)\to S_r(x_0,X)\times(r,t]$
given by $\psi(x)=(c_{s,r}(x),s)$, where $s=d(x,x_0)$. Note that,
due to geodesic extension property in $A_{r,t}$, the map $\psi$ is surjective.
Since $\psi$ extends, with the same formula, to a continuous map
$\bar\psi:A_{r,t}\to S_r(x_0,X)\times[r,t]$ between compact spaces,
and since $\bar\psi^{-1}(S_r(x_0,X)\times\{t\})=S_t(x_0,X)$,
it is easy to deduce that $\psi$ is proper.
Moreover, for each $(y,s)\in S_r(x_0,X)\times(r,t]$ the inverse image
$\psi^{-1}((y,s))$ coincides with the set $c_{s,r}^{-1}(y)$,
and hence $(C_n)$ implies that this set is cell-like. Consequently, $\psi$
is a cell-like map of manifolds. Consider the function 
$\delta:S_r(x_0,X)\times(r,t]\to R_+$ given by $\delta(y,s)=s-r$.
By Approximation Theorem, there is a homeomorphism
$h:A_{r,t}\setminus S_r(x_0,X)\to S_r(x_0,X)\times(r,t]$ such that
$d_X(h(y,s),\psi(y,s))<\delta(y,s)=s-r$ for each $(y,s)\in S_r(x_0,X)\times(r,t]$.
Consequently, if an argument $z\in A_{r,t}\setminus S_r(x_0,X)$ converges to
some $z_0\in S_r(x_0,X)$, then $h(z)$ converges to $(z_0,r)$.
Hence $h$ can be extended to a continous map $H:A_{r,t}\to
S_r(x_0,X)\times[r,t]$ by putting $H(x)=(x,r)$ for $x\in S_r(x_0,X)$.
Moreover, $H$ is easily seen to be a bijection, and since its domain is compact,
it is a homeomorphism. Finally, since homeomorphisms preserve the boundary,
$h$ maps $S_t(x_0,X)$ onto the set $S_r(x_0,X)\times\{t\}$, and the same
is true for $H$. It is also clear that
$H$ maps $S_r(x_0,X)$ onto the set
$S_r(x_0,X)\times\{r\}$, and this completes the proof.

\medskip
\noindent
{\bf Fact 6.}
{\it Under assumptions as in Lemma B, suppose additionally that
for each $x\in X$ such that $d(x,x_0)=r$
the metric link
$\hbox{\rm Lk}(x,X)$ (viewed as a piecewise spherical complex) is
PL homeomorphic to the standard PL sphere $S^{n-1}$. 
Then the sphere $S_r(x_0,X)$ is a closed $(n-1)$-manifold which is
bi-collared in $X$ (i.e. $S_r(x_0,X)$ has a neighbourhood
homeomorphic to the product $Y\times(-1,1)$, with $Y\times\{0\}$
corresponding to $S_r(x_0,X)$).}

\medskip\noindent
{\bf Proof:}
The assumption concerning links implies that there is $\epsilon>0$
such that for any $x\in X$ satisfying $d(x,X)\in(r-\epsilon,r+\epsilon)$
the link
$\hbox{\rm Lk}(x,X)$ is
PL homeomorphic to the standard PL sphere $S^{n-1}$. 
This follows by observing that the sphere $S_r(x_0,X)$ is compact
and each of its points has a cone neighbourhood in $X$.
Then, it follows from $(D_n)$ that the annular region $A_{r-\epsilon/2,r+\epsilon/2}$ is a bi-collared neighbourhood 
of $S_r(x_0,X)$ in $X$, as required.

\bigskip
We now prove part (1) of $(A_n)$.
By existence of a cone neighbourhood,
there is $\epsilon>0$ such that the closed ball $B_\epsilon(v,L)$
is homeomorphic to the $n$-disk $B_n$, with the sphere $S_\epsilon(v,L)$
corresponding to the boundary. By Fact 5, the annular region
$A_{\epsilon,r}=B_r(v,L)\setminus B_\epsilon^\bullet(v,L)$
is homeomorphic to $S_\epsilon(v,L)\times[\epsilon,r]\cong 
S^{n-1}\times[\epsilon,r]$, with $S_\epsilon(v,L)$ coinciding with
one of the boundary components of this region.
Consequently, the pair $(B_r(v,L),S_r(v,L))$ is homeomorphic to 
the pair $(B^n,\partial B^n)$, as required. It is also an obvious consequence
of Fact 6 that $B_r(v,L)$ is collared in $L$. This completes the proof of part (1).

\medskip
We now turn to proving part (2) of $(A_n)$.
Note that for $r\in(0,\pi/2)$ the first assertion of part (2) follows directly from
part (1). Thus, we only need to justify this assertion for $r\in[\pi/2,\pi]$.
Given any such $r$, consider any homeomorphism $\lambda:[0,r)\to[0,\pi/4)$
such that $\lambda(t)<t$ for all $t>0$ in the domain. 
Define the map $f:B^\bullet_r(x_0,X)\to B^\bullet_{\pi/4}(x_0,X)$
by letting the restriction of $f$ to any sphere $S_t(x_0,X)$, $t<r$, be the
geodesic projection $c_{t,\lambda(t)}$ which maps $S_t(x_0,X)$ onto 
$S_{\lambda(t)}(x_0,X)$.
It follows then from $(C_n)$ that $f$ is a cell-like map between $n$-manifolds.
Hence, by Approximation Theorem, $f$ can be approximated by homeomorphisms.
In particular, the open ball $B^\bullet_r(x_0,X)$ is homeomorphic to
$B^\bullet_{\pi/4}(x_0,X)$, and thus also to the open $n$-disk.

The second assertion of part (2) follows immediately from the first one. 

\medskip
To prove part (3) of $(A_n)$, put $C=L\setminus B^\bullet_r(v,L)$.
Note that, due to part (2) of $(A_n)$, the complement of $C$ in the PL $n$-sphere
$L$ is an open $n$-disk. This property of $C$ is known under the term 
{\it point-like}. It is also known that any point-like subset of a sphere is cellular
(see e.g. [Ed], Theorem on page 114), so the assertion follows
simply by referring to this result.

%\vfill\break

\bigskip
\noindent
\centerline{\bf References}

\medskip

\itemitem{[AS]} F. Ancel, L. Siebenmann,
{\it The construction of homogeneous homology manifolds},
Abstracts Amer. Math. Soc. 6 (1985), 92.

\itemitem{[Ar]} S. Armentrout, {\it Cellular decompositions of 3-manifolds
that yield 3-manifolds}, Mem. Amer. Math. Soc. 107 (1971).

\itemitem{[BH]} M. Bridson, A. Haefliger, 
Metric Spaces of Non-Positive
Curvature, Grundlehren der mathematischen Wissenschaften 319, Springer,
1999.

\itemitem{[Br]} M. Brown, {\it Some applications of an appropximation
theorem for inverse limits}, Proc. Amer. Math. Soc. 11 (1960), 478--481. 

\itemitem{[C]} R. Coulon, {\it Asphericity and small 
cancellation theory for rotation families of groups},
Groups Geom. Dyn. 5 (2011), 729--765.

\itemitem{[CD]} R. Charney, M. Davis,  {\it Strict hyperbolization}, 
Topology 34 (1995), 329--350.

\itemitem{[Da]} M. Davis, {The geometry and topology 
of Coxeter groups}, {London Mathematical Society Monographs Series},
{vol. 32}, {Princeton University Press}, {Princeton}, {2008}.

\itemitem{[DJ]} M. Davis, T. Januszkiewicz, 
{\it Hyperbolization of polyhedra}, J. Differential Geometry 34 (1991),
347--388.

\itemitem{[Ed]} R. D. Edwards, {\it The topology of manifolds
and cell-like maps}, Proc. Int. Congress of Math. Helsinki, 1978.  

\itemitem{[En]} R. Engelking, General Topology, Heldermann Verlag, 
Berlin, 1989.

\itemitem{[F]} H. Fischer, {\it Boundaries of right--angled Coxeter 
groups with manifold nerves},
Topology 42 (2003), 423--446.

\itemitem{[FG]} H. Fischer, C. Guilbault, {\it On the fundamental groups
of trees of manifolds}, Pacific J. Math. 221 (2005), 49--79. 

\itemitem{[FM]} K. Fujiwara, J. F. Manning, {\it CAT$(0)$ and CAT$(-1)$
fillings of hyperbolic manifolds,}
J. Differential Geom.  85 (2010), 229--270. 

\itemitem{[J]} W. Jakobsche, {\it Homogeneous cohomology manifolds which are inverse limits},
Fundamenta Mathematicae 137 (1991), 81--95.

\itemitem{[La]} R. C. Lacher, {\it Cell-like mappings, I},
Pacific J. Math. 30 (1969), 717--731.

\itemitem{[Mo]} R. L. Moore, {\it Concerning upper semi-continuous
collections of continua}, Trans. Amer. Math. Soc. 27 (1925), 416--428.

\itemitem{[MS]} L. Mosher, M. Sageev, {\it Non-manifold hyperbolic
groups of high cohomological dimension,} preprint 1997.

\itemitem{[PS]} P. Przytycki, J. \'Swi\c atkowski, 
{\it Flag-no-square triangulations and Gromov boundaries
in dimension 3}, Groups, Geometry \& Dynamics  3 (2009), 453--468.

\itemitem{[Qu]} F. Quinn, {\it Ends of maps. III: dimensions 4 and 5},
J. Diff. Geom. 17 (1982), 503--521.

\itemitem{[Si]} L. Siebenmann, {\it Approximating cellular maps
by homeomorphism}, Topology 11 (1973), 271--294.

\itemitem{[St]} P.R. Stallings, {\it An extension of Jakobsches construction
of $n$--ho\-mo\-geneous continua to the
nonorientable case}, in Continua (with the Houston Problem Book),
ed. H. Cook, W.T. Ingram, K. Kuperberg, A. Lelek, P. Minc,
Lect. Notes in Pure and Appl. Math. vol. 170 (1995), 347--361.

\itemitem{[Sw]} J. \'Swi\c atkowski, {\it Trees of metric compacta
and trees of manifolds},
preprint 2013.

\itemitem{[Z]} P. Zawi\'slak, {\it Trees of manifolds and boundaries of 
systolic groups}, Fund. Math. 207 (2010), 71--99.

\bigskip

\noindent  Instytut Matematyczny, Uniwersytet Wroc\l awski, 

\noindent pl. Grunwaldzki 2/4, 50-384 Wroc\l aw, Poland

\smallskip
\noindent E-mail: {\tt Jacek.Swiatkowski@math.uni.wroc.pl}

\bye